\documentclass[preprint,12pt]{elsarticle}

\makeatletter
\def\ps@pprintTitle{%
	\let\@oddhead\@empty
	\let\@evenhead\@empty
	\let\@oddfoot\@empty
	\let\@evenfoot\@empty
}
\renewcommand\section{\@startsection{section}{1}{\z@}%
  {-18\p@ \@plus -6\p@ \@minus -3\p@}%
  {9\p@ \@plus 6\p@ \@minus 3\p@}%
  {\normalfont\normalsize\bfseries}}
\renewcommand\subsection{\@startsection{subsection}{2}{\z@}%
  {-12\p@ \@plus -6\p@ \@minus -3\p@}%
  {6\p@ \@plus 6\p@ \@minus 3\p@}%
  {\normalfont\normalsize\bfseries}}
\makeatother

\usepackage[table]{xcolor}
\usepackage{graphicx}
\usepackage{tikz}
\usetikzlibrary{backgrounds,arrows}
\usepackage{pgfplots}
\pgfplotsset{compat=1.15}
\usepackage{mathrsfs}
\usepackage[top=2cm,bottom=2cm,left=1.5cm,right=1.5cm,headheight=14pt,headsep=0.25in,footskip=0.5in]{geometry}
\usepackage{float}
\usepackage{amsmath,amssymb,amsthm}
\usepackage{tabularx}
\usepackage{comment}
\usepackage{etoolbox}
\usepackage[colorlinks,linkcolor=blue]{hyperref}
\usepackage[nameinlink, capitalize]{cleveref}

\def\arraystretch{1.1}
\def\m#1#2#3{\raise 0.2ex\hbox{
    ${#1_{\displaystyle\textcolor{black} #2}(\textcolor{red}{#3})}$}}

\definecolor{pl0}{HTML}{E6194B}
\definecolor{pl1}{HTML}{F58231}
\definecolor{pl2}{HTML}{FFD900}
\definecolor{pl3}{HTML}{3CB44B}
\definecolor{pl4}{HTML}{2EC4E6}
\definecolor{pl5}{HTML}{2F4BD0}
\definecolor{pl6}{HTML}{911EB4}
\definecolor{pl7}{HTML}{E500E5}
\definecolor{pl8}{HTML}{808080}
\definecolor{skyblue}{HTML}{A8D8F0}
\definecolor{softred}{HTML}{F4A6A0}

\tikzset{
  vtx/.style={circle,draw,line width=0.8pt,inner sep=0.5pt,minimum size=0.55cm,font=\tiny,fill=white},
  redge/.style={blue!70!black,line width=1.5pt},
  cedge/.style={red!70!black,line width=1.5pt},
  sedge/.style={black,line width=1.5pt,dash pattern=on 4pt off 2.5pt},
  uvtx/.style={circle,draw,line width=0.8pt,inner sep=0.5pt,minimum size=0.6cm,font=\scriptsize},
  rowe/.style={black,line width=1.5pt},
  cole/.style={black,line width=1.5pt,dash pattern=on 5pt off 3pt}
}

\newcommand{\qcube}[8]{%
\begin{tikzpicture}[scale=1.1]
  \node[vtx] (o1) at (0,2.2) {#1};
  \node[vtx] (o2) at (2.2,2.2) {#2};
  \node[vtx] (o3) at (2.2,0) {#3};
  \node[vtx] (o4) at (0,0) {#4};
  \node[vtx] (i1) at (0.65,1.55) {#5};
  \node[vtx] (i2) at (1.55,1.55) {#6};
  \node[vtx] (i3) at (1.55,0.65) {#7};
  \node[vtx] (i4) at (0.65,0.65) {#8};
  \begin{scope}[on background layer]
    \draw[redge] (o1)--(o2) (o4)--(o3) (i1)--(i2) (i4)--(i3);
    \draw[sedge] (o1)--(o4) (o2)--(o3) (i1)--(i4) (i2)--(i3);
    \draw[cedge] (o1)--(i1) (o2)--(i2) (o3)--(i3) (o4)--(i4);
  \end{scope}
  \node[font=\small] at (1.1,-0.7) {$Q_3$};
\end{tikzpicture}}

\newtheorem{thm}{Theorem}[section]
\newtheorem{prop}[thm]{Proposition}
\newtheorem{lem}[thm]{Lemma}
\newtheorem{cor}[thm]{Corollary}
\newtheorem{defin}[thm]{Definition}
\newtheorem{rem}[thm]{Remark}

\newtheorem{conj}[thm]{Conjecture}

\makeatletter
\newcommand\noindentafter[1]{\AfterEndEnvironment{#1}{\@afterindentfalse\@afterheading}}
\makeatother
\noindentafter{thm}\noindentafter{prop}\noindentafter{lem}\noindentafter{cor}
\noindentafter{defin}\noindentafter{rem}\noindentafter{con}\noindentafter{ex}
\noindentafter{ob}\noindentafter{nota}\noindentafter{conj}\noindentafter{proof}

\journal{Journal of Combinatorial Theory, Series A}
\begin{document}

\begin{frontmatter}

\title{On the Chromatic Number of Cayley Tables of Dihedral Groups}

\author[1]{Kianoosh Shokri} 

\affiliation[1]{organization={Department of Mathematics and Statistics, University of Ottawa},
            city={Ottawa},
            country={Canada}}

\author[2]{Masoomeh Akbari } 

\affiliation[2]{organization={Department of Mathematics and Statistics, University of British Columbia},
            city={Vancouver},
           country={Canada}}

\begin{abstract}
A coloring of a Latin square $L$ is an assignment of $k$ colors to its cells such that no two cells sharing a row, a column, or a symbol receive the same color. The smallest $k$ for which such a coloring exists is called the \emph{chromatic number} of $L$, denoted by $\chi(L)$. The Cayley table of a finite group $G$ of order $n$ is a Latin square, which we denote by $L_G$. It is known that  $\chi(L_G)\ne n+1$, so either $\chi(L_G)=n$ or $\chi(L_G)\ge n+2$. Moreover, $\chi(L_G)=n$ if and only if every Sylow $2$-subgroup of $G$ is trivial or non-cyclic, which completely characterizes the finite groups whose Cayley tables attain the lower bound $n$.

In this paper, we consider the dihedral group $D_n$ of order $2n$, where $n$ is odd. Since $D_n$ has a cyclic Sylow $2$-subgroup, it follows that $\chi(L_{D_n})\ge 2n+2$. We prove that this bound is sharp, that is, $\chi(L_{D_n})=2n+2$.
\end{abstract}

\begin{keyword}
Latin square \sep chromatic number \sep Sylow subgroups  \sep dihedral groups \sep Cayley table

\end{keyword}

\end{frontmatter}


\section{Introduction}

A \emph{Latin square} $L = [L_{r,c}]$ of \emph{order} $n$ is an $n \times n$ array on $n$ symbols in which each symbol occurs exactly once in each row and each column. If $L_{r,c} = s$, we denote the cell in row $r$ and column $c$ by the triple $(r, c; s)$. Two distinct cells $(r, c; s)$ and $(r', c'; s')$ are \emph{independent} if $r \neq r'$, $c \neq c'$, and $s \neq s'$; otherwise, they are \emph{dependent}.

A \emph{partial transversal} of \emph{length} $k$ is a set of $k$ pairwise independent cells, and a \emph{transversal} is a partial transversal of length $n$. More generally, a \emph{$k$-plex} is a set of $kn$ cells containing exactly $k$ cells from each row, each column, and each symbol; thus a transversal is a $1$-plex.

The \emph{Latin square graph} $\Gamma(L)$ is the simple graph whose vertices are the cells of $L$, where two distinct cells are adjacent if they are dependent. It is well known that $\Gamma(L)$ is strongly regular with parameters $(n^2, 3(n-1), n, 6)$; that is, $\Gamma(L)$ has $n^2$ vertices, each of degree $3(n-1)$, any two adjacent vertices have exactly $n$ common neighbors, and any two nonadjacent vertices have exactly $6$ common neighbors.

The independent sets of $\Gamma(L)$ are precisely the partial transversals of $L$. Hence a proper coloring of $\Gamma(L)$ is a partition of the cells of $L$ into partial transversals, and the chromatic number of $\Gamma(L)$ is the minimum number of partial transversals needed to partition $L$. This number is called the \emph{chromatic number} of $L$ and is denoted by $\chi(L)$.

Since every partial transversal has length at most $n$, at least $n$ partial transversals are needed to partition the $n^2$ cells of $L$; hence $\chi(L) \geq n$. Moreover, $\chi(L) = n$ if and only if the cells of $L$ can be partitioned into $n$ transversals, which is equivalent to $L$ having an orthogonal mate~\cite{mahmoodian2016}. Here, a Latin square $M$ of order $n$ is an \emph{orthogonal mate} of $L$ if the $n^2$ ordered pairs $(L_{r,c}, M_{r,c})$ are all distinct.

Besharati et al.~\cite{mahmoodian2016} proposed the following conjecture and verified it for all Latin squares of order $n \leq 8$.

\begin{conj}[{\cite{mahmoodian2016}}] \label{conj1}
    Every Latin square $L$ of order $n$ satisfies
    \begin{equation} \label{eq1}
        \chi(L) \leq \begin{cases}
        n + 1 & \text{if $n$ is odd,} \\
        n + 2 & \text{if $n$ is even.}
        \end{cases}
    \end{equation}
\end{conj}


\begin{prop}[{\cite{mahmoodian2016}}] \label{prop1}
    Let $L$ be a Latin square of order $n$. If $\chi(L) \leq n + 1$, then $L$ has a transversal. If $\chi(L) \leq n + 2$, then $L$ has a partial transversal of length at least $n-1$.
\end{prop}

\begin{proof}
    Fix a proper coloring of $L$ with $\chi(L)$ colors, and let $C_1, \dots, C_{\chi(L)}$ be its color classes. Each $C_i$ is a partial transversal, and these classes partition the $n^2$ cells of $L$.

    Suppose $\chi(L) \leq n+1$ but $L$ has no transversal. Then $|C_i| \leq n-1$ for every $i$, so
    $
        n^2 = \sum_{i=1}^{\chi(L)} |C_i| \leq (n+1)(n-1) = n^2 - 1,
    $
    a contradiction.

    Suppose $\chi(L) \leq n+2$ but $L$ has no partial transversal of length at least $n-1$. Then $|C_i| \leq n-2$ for every $i$, so
    $
        n^2 = \sum_{i=1}^{\chi(L)} |C_i| \leq (n+2)(n-2) = n^2 - 4,
    $
    a contradiction.
\end{proof}


As discussed in~\cite{mahmoodian2016}, if \cref{conj1} is true, then by \cref{prop1} it implies two long-standing conjectures on transversals, the Brualdi--Stein conjecture and Ryser's conjecture, both stated below.

\begin{conj}[{\cite{denes1974,stein1975}}] \label{conj:BS}
    Every Latin square of even order $n$ contains a partial transversal of length $n-1$. 
\end{conj}
\begin{conj}[{\cite{ryser1967}}] \label{conj:Ryser}
    Every Latin square of odd order contains a transversal.
\end{conj}
Both conjectures remain open in general, although Montgomery~\cite{montgomery2023} recently proved \cref{conj:BS} for all sufficiently large $n$.


\medskip

A special type of Latin square is obtained by using a finite group. Let $(G, \circ)$ be a finite group of order $n$ with a fixed enumeration $G = \{g_0, g_1, \ldots, g_{n-1}\}$. The \emph{Cayley table} of $G$ is the $n \times n$ array, denoted by $L_G = [a_{i,j}]$, where $a_{i,j} = g_i \circ g_j$. It is easy to see that $L_G$ is a Latin square.  Note that a different enumeration of $G$ only permutes the rows and columns of $L_G$, so $\chi(L_G)$ does not depend on the choice of enumeration.

The structure of a finite group makes the study of $\chi(L_G)$ easier and helps toward proving \cref{conj1} for special cases of Latin squares of finite groups. Although proving \cref{conj1} in general is still difficult for Latin squares of groups, significant improvements have been obtained in characterizing groups with respect to their chromatic numbers.  For a finite group $G$ of order $n$, \cref{thm8} implies that $\chi(L_G) = n$ or $\chi(L_G) \geq n + 2$. In addition, \cref{thm3} provides a necessary and sufficient condition to determine when $\chi(L_G) = n$. 

To make progress toward a proof of \cref{conj1} in the special case of the Latin square of a group $G$ of order $n$, we aim to show that $\chi(L_G) = n$ or $\chi(L_G) = n + 2$. The main goal of this paper is to study the chromatic number of the Latin square of the dihedral group $D_n$. 

Studying the chromatic number of Latin squares of groups requires tools from algebra. In \cref{pre}, we review basic concepts and theorems from algebra and demonstrate their connection with the coloring of Latin squares of groups. In \cref{exist}, we study certain groups $G$ of order $n$ with $\chi(L_G) = n$. In particular, in \cref{cor4} we prove that if $n$ is even, then $\chi(L_{D_n}) = 2n$. In \cref{preresult}, we review the chromatic numbers of Latin squares of some other finite groups with cyclic Sylow 2-subgroups. In particular, we see that $\chi(L_G) = n + 2$ for any commutative group $G$ of order $n$ with a cyclic Sylow 2-subgroup. In \cref{main}, we study the chromatic number of the Latin square of the dihedral group $D_n$ for odd $n$. \Cref{thmmain}, the main result of this paper, states that if $n$ is odd, then $\chi(L_{D_n}) = 2n + 2$. In \cref{future}, we summarize our results and discuss directions for future work.

\section{Preliminaries} \label{pre}

In this section, we provide the necessary background on algebra and demonstrate its connection with the coloring of Latin squares of groups. The goal of this section is to review a fundamental theorem (\cref{thm3}), which determines whether the Cayley table of a finite group has a transversal based on its Sylow $2$-subgroups.

\subsection{Background on algebra}

Let $G$ be a finite group with $n$ elements. Then the order or size of $G$ is denoted by $|G| = n$. Let $a \in G$ and let $e$ be the identity element of $G$. The order of $a$ is the smallest positive integer $k$ such that $a^k = e$ and is denoted by $|a| = k$.
Let $p$ be a prime number. A group $G$ is a \emph{$p$-group} if the order of every element of $G$ is a power of $p$. A subgroup of $G$ is a \emph{$p$-subgroup} if it is a $p$-group.
A \emph{Sylow $p$-subgroup} $P$ of a group $G$ is a maximal $p$-subgroup of $G$.
A subgroup $H$ of $G$ is \emph{normal} in $G$ if $gH = Hg$ for all $g \in G$. 
 A nontrivial group $G$ whose only normal subgroups are $\{e\}$ and $G$ is called a \emph{simple group}.

\begin{thm}[{\cite{jud}}] \label{thmlag}
    Let $G$ be a finite group and $H$ a subgroup of $G$. Then $|H|$ divides $|G|$.
\end{thm}

\begin{thm} [{\cite{jud}}] \label{thm4}
    Let $G$ be a finite group and $p$ a prime number that divides the order of $G$. Then $G$ contains a subgroup of order $p$. 
\end{thm}

\begin{thm}[{\cite{dummit2003}}] \label{thm7}
    Let $G$ and $H$ be two finite groups with Sylow $p$-subgroups $P$ and $Q$, respectively. Then $P \times Q$ is a Sylow $p$-subgroup of $G \times H$. 
\end{thm}

\begin{thm}[{\cite{scott1987}}] \label{thm6}
    Let $G$ be a finite group with $|G| = 2^km$, $\gcd(2, m) = 1$. If $G$ has a cyclic Sylow $2$-subgroup, then $G$ has a unique normal subgroup of order $m$. 
\end{thm}

\begin{rem} \label{sylowsize}
    Let $G$ be a group with $|G| = 2^km$, where $\gcd(2, m) = 1$. If $G$ has a Sylow $2$-subgroup $H$, then $|H| = 2^k$. 
\end{rem}

\subsection{Latin squares of groups and complete mapping}

In this section, we see how the algebraic approach helps us determine the chromatic number of $L_G$ for a finite group $G$. 

\begin{thm}[{\cite{wanless2011}}] \label{thm1}
    If the Cayley table of a finite group has a single transversal, then it has a decomposition into disjoint transversals. 
\end{thm}

\cref{thm1} shows that if $L_G$ has one transversal, then $\chi(L) = n$. This leads to the characterization of finite groups with and without a transversal.

\begin{defin}
    Let $(G, \circ)$ be a finite group. A permutation $\theta$ of $G$ is a \emph{complete mapping} if the map $\eta : G \to G$ defined by $\eta(x) = x \circ \theta(x)$ is also a permutation of $G$.
\end{defin}

A complete mapping is an equivalent concept of a transversal in groups.

\begin{thm}[{\cite{wanless2011}}] \label{thm2}
    Let $(G, \circ)$ be a finite group, $L_G$ its Cayley table, and $\theta : G \longrightarrow G$ a permutation of $G$. Let $\eta : G \longrightarrow G$ be defined by $\eta(x) = x \circ \theta(x)$. Then $\theta$ is a complete mapping if and only if selecting, in each row $x$, the entry $x \circ \theta(x) = \eta(x)$ in column $\theta(x)$ gives a transversal of $L_G$.
\end{thm}

Studying transversals of $L_G$ is equivalent to studying complete mappings of $G$. \Cref{thm5} provides a necessary and sufficient condition for a finite commutative group $G$ to have a complete mapping, and \cref{thm3} extends this to all finite groups by giving several equivalent conditions for the existence of a complete mapping.

\begin{thm} [{\cite{Paige1947}}] \label{thm5}
    Let $G$ be a finite commutative group. Then $G$ does not have a complete mapping if and only if $G$ has a unique element of order $2$.
\end{thm}

\begin{thm}[{\cite{wanless2011}}] \label{thm3}
    Let $(G, \circ)$ be a finite group and $L_G$ its Cayley table. Then the following conditions are equivalent. 
    \begin{enumerate}
        \item $L_G$ has a transversal. 
        \item $L_G$ can be decomposed into disjoint traversals. 
        \item There exists a Latin square orthogonal to $L_G$. 
        \item There is some ordering of the elements of $G$, say $g_0, g_1, \ldots, g_{n-1}$, such that $g_0 \circ g_1 \circ \ldots \circ g_{n-1} = e$, where $e$ denotes the identity element of $G$. 
        \item The Sylow 2-subgroups of $G$ are trivial or non-cyclic. 
    \end{enumerate}
\end{thm}

\Cref{thm3} gives a powerful tool to determine whether $L_G$ has a transversal by studying the Sylow $2$-subgroups of $G$. In addition, if $G$ is commutative, \cref{thm5} is often easier to apply for the same purpose.

\begin{rem} \label{rem1}
    Let $G$ be an arbitrary finite group with more than one non-trivial Sylow $2$-subgroup. Since all Sylow $2$-subgroups are conjugate, one Sylow $2$-subgroup is cyclic if and only if all Sylow $2$-subgroups are cyclic. 
\end{rem}

By \cref{thm3}, $L_G$ has a transversal if and only if all Sylow $2$-subgroups are trivial or non-cyclic. However, by \cref{rem1}, it suffices to study only one Sylow $2$-subgroup. 

\section{Existence of a complete mapping for certain finite groups} \label{exist}

In this section, we study the existence of a complete mapping for certain finite groups. The main goal of this section is to review which groups have a transversal and which groups do not. Some of the results that follow were originally obtained by various, often lengthy, methods before \cref{thm3} was established. Here, we present short proofs of these results as direct applications of \cref{thm3}. If $L_G$ has a transversal, then $\chi(L_G) = n$, while if $L_G$ has no transversal, then $\chi(L_G) > n$. Hence, characterizing finite groups based on the existence of a transversal immediately gives us information about the chromatic number.

The first general result on the existence or non-existence of a complete mapping is given in \cref{cor1}. It was originally proved by Mann~\cite{Mann1944} and Paige~\cite{Paige1951} using direct methods, whereas here it follows easily from \cref{thm3}.

\begin{cor} \label{cor1}
    Let $G$ be a group of order $n$. 
    \begin{enumerate}
        \item If $n$ is odd, then $G$ has a complete mapping. 
        \item If $n = 4k + 2$, for $k \geq 0$, then $G$ does not have a complete mapping. 
    \end{enumerate}
\end{cor}

\begin{proof}
    Let $|G| = n$ be odd; then the Sylow $2$-subgroup of $G$ is trivial. By \cref{thm3}, $G$ has a complete mapping. 
    
    For the second part, let $|G| = 2(2k+1)$. By \cref{thm4}, $G$ has a subgroup of order $2$ isomorphic to $\mathbb{Z}_2$. Since $\gcd(2, 2k + 1) = 1$, this subgroup is maximal and hence a cyclic Sylow $2$-subgroup. Therefore, by \cref{thm3}, $G$ does not have a complete mapping.  
\end{proof}

\Cref{cor2} gives results on the existence of complete mappings for certain commutative groups.

\begin{cor} \label{cor2}
The following statements hold: 
\begin{enumerate}[(i)]
    \item The group $\mathbb{Z}_{2k}$ does not have a complete mapping. \label{cor21}
    \item The group $\mathbb{Z}_n \times \mathbb{Z}_2$ has a complete mapping if and only if $n = 2k$. 
    \item Let $G$ be a finite commutative group with a subgroup isomorphic to $\mathbb{Z}_2^k$ for some $k \geq 2$. Then $G$ has a complete mapping. 
\end{enumerate}
    
\end{cor}

\begin{proof}
   (i) The unique element of order $2$ in $\mathbb{Z}_{2k}$ is $k$. Hence, by \cref{thm5}, $\mathbb{Z}_{2k}$ does not have a complete mapping.

   (ii) For the sufficient condition, let $n = 2k$. Then $(k, 0)$ and $(k, 1)$ are two distinct elements of order $2$, and by \cref{thm5}, $\mathbb{Z}_{2k} \times \mathbb{Z}_2$ has a complete mapping. For the necessary condition, let $n = 2k + 1$. Then $(0, 1)$ is the unique element of order $2$, and by \cref{thm5}, $\mathbb{Z}_{2k+1} \times \mathbb{Z}_2$ does not have a complete mapping.

   (iii) Since $G$ has a subgroup isomorphic to $\mathbb{Z}_2^k$, and $(1, 0, \ldots, 0)$ and $(0, 1, 0, \ldots, 0)$ are two distinct elements of order $2$ in $\mathbb{Z}_2^k$, the group $G$ has at least two elements of order $2$. Therefore, by \cref{thm5}, $G$ has a complete mapping.
\end{proof}

In the following, we study certain non-commutative groups regarding the existence of a complete mapping.  

\begin{cor} \label{cor3}
    The symmetric group $S_n$ for $n \geq 4$ and the alternating group $A_n$ for $n \in \mathbb{N}$ have a complete mapping.
\end{cor}

\begin{proof}

    First, we prove that $S_n$, $n \geq 5$, has a complete mapping. The group $S_4$ is discussed later. Note that the unique normal subgroup of $S_n$, $n \geq 5$ is $A_n$. Let $|S_n| = 2^km$, $\gcd(2, m) = 1$. Then $|A_n| = 2^{k-1}m$. Hence, $S_n$ does not have a normal subgroup of order $m$. So, by \cref{thm6}, $S_n$ does not have a cyclic Sylow $2$-subgroup. Thus, Sylow $2$-subgroups of $S_n$ must be trivial or non-cyclic. Therefore, by \cref{thm3}, $S_n$, $n \geq 5$ has a complete mapping. 
    The group $S_4$ has two proper nontrivial normal subgroups, $A_4$ and $\{ e, (1 2)(3 4), (1 3)(2 4), (1 4)(2 3) \}$. Since $|S_4| = 2^3 \cdot 3$, we have $m = 3$, and $S_4$ does not have a normal subgroup of order $3$. So $S_4$ does not have a cyclic Sylow $2$-subgroup. Therefore, $S_4$ has a complete mapping.

Next, we prove that $A_n$ for $n \in \mathbb{N}$ has a complete mapping. For $n \leq 3$, $|A_n|$ is odd, so the Sylow $2$-subgroups of $A_n$ are trivial. Therefore, by \cref{thm3}, $A_n$ has a complete mapping. For $n \geq 5$, $A_n$ is a simple group. Let $|A_n| = 2^km$, $\gcd(2, m) = 1$. Since $1 < m < |A_n|$, $A_n$ does not have a normal subgroup of order $m$, which by \cref{thm6} implies that $A_n$ does not have a cyclic Sylow $2$-subgroup. Therefore, by \cref{thm3}, $A_n$ for $n \geq 5$ has a complete mapping. The unique proper nontrivial normal subgroup of $A_4$ is $\{ e, (1 2)(3 4), (1 3)(2 4), (1 4)(2 3) \}$. Since $|A_4| = 2^2 \cdot 3$, we have $m = 3$, and $A_4$ does not have a normal subgroup of order $3$. So, by \cref{thm6}, $A_4$ does not have a cyclic Sylow $2$-subgroup. Therefore, by \cref{thm3}, $A_4$ has a complete mapping.
 
\end{proof}

Let $D_n = \left < r, s : r^{n} = e, s^2 = e, (rs)^2 = e \right >$ be the dihedral group of order $2n$. The elements of $D_n$ are $\{r^k, r^ks : 0 \leq k \leq n-1 \}$ with the following relation: 
\begin{equation*}
    r^ir^j = r^{i + j}, \quad r^isr^js = r^{i-j}, 
\end{equation*}
where $0 \leq i, j \leq n-1$ and all exponents are modulo $n$. 
\begin{cor} \label{cor4}
    The dihedral group $D_n$ has a complete mapping if and only if $n$ is even.  
\end{cor}

\begin{proof}
    First, we prove that if $n$ is even, then $D_n$ has a complete mapping. Let $n = 2k$ and $|D_{2k}| = 4k = 2^lm$, where $\gcd(2, m) = 1$ and $l \geq 2$. Suppose, for contradiction, that $H$ is a cyclic Sylow $2$-subgroup of $D_{2k}$. By \cref{sylowsize}, $|H| = 2^l$. Since $H$ is cyclic, $H = \langle r^is \rangle$ or $H = \langle r^i \rangle$ for some $0 \leq i \leq 2k-1$. Since $|\langle r^is \rangle| = 2 < 2^l$, we have $H = \langle r^i \rangle$ with $0 < i \leq 2k-1$. Let $d = \gcd(i, 2k)$. Then $\langle r^i \rangle = \langle r^d \rangle$, and since $d \mid 2k$, we have $2k = dj$ for some positive integer $j$. So, $H = \{ e, r^d, r^{2d}, \ldots, r^{(j-1)d} \}$, and hence $j = |H| = 2^l$. Thus, $4k = 2dj = 2^{l+1}d$. Since $4k = 2^lm$, we obtain $m = 2d$, which contradicts $\gcd(2, m) = 1$. Therefore, $D_{2k}$ has no cyclic Sylow $2$-subgroup. Since $l \geq 2$, the Sylow $2$-subgroups of $D_{2k}$ are nontrivial, and hence non-cyclic. Therefore, by \cref{thm3}, $D_{2k}$ has a complete mapping.

    For the reverse, suppose $n$ is odd. Then $|D_n| = 2n$, where $2$ is the largest power of $2$ dividing $2n$. Consider $\langle s \rangle = \{ e, s \}$, which is a cyclic subgroup of order $2$. By \cref{thmlag}, every $2$-subgroup of $D_n$ has order at most $2$, so $\langle s \rangle$ is a maximal $2$-subgroup, that is, a cyclic Sylow $2$-subgroup of $D_n$. Therefore, by \cref{thm3}, $D_n$ does not have a complete mapping.
\end{proof}

\begin{cor} \label{cor44}
    The group $D_n \times \mathbb{Z}_{2m + 1}$ has a complete mapping if and only if $n$ is even. 
\end{cor}

\begin{proof}
    First, we prove that if $n$ is even, then $D_n \times \mathbb{Z}_{2m + 1}$ has a complete mapping. When $n$ is even, by \cref{cor4}, $D_n$ has a complete mapping and does not have a non-trivial cyclic Sylow $2$-subgroup. Thus, $D_n \times \mathbb{Z}_{2m + 1}$ does not have a non-trivial cyclic Sylow $2$-subgroup. Therfore, $D_n \times \mathbb{Z}_{2m + 1}$ has a complete mapping. 
    For the reverse, suppose $n$ is odd. Then $D_n$ has Sylow $2$-subgroup isomorphic to $\mathbb{Z}_2$. Since Sylow $2$-subgroup of $\mathbb{Z}_{2m+1}$ is trivial, by \cref{thm7}, Sylow $2$-subgroup of $D_n \times \mathbb{Z}_{2m + 1}$ is isomorphic to $\mathbb{Z}_2$ and hence has a cyclic Sylow $2$-subgroup which by \cref{thm3} does not have a complete mapping. 
\end{proof}

Let $Q_m = \langle a, b : a^{2m} = e, b^2 = a^{m}, ba = a^{-1}b \rangle$, $m \geq 2$, be the generalized quaternion group of order $4m$. The elements of $Q_m$ are $\{a^k, a^kb : 0 \leq k \leq 2m-1 \}$ with the following relations:
\begin{equation*}
    a^ia^j = a^{i + j}, \quad a^i a^jb = a^{i+j}b, \quad a^iba^j = a^{i-j}b, \quad a^iba^jb = a^{i-j+m},
\end{equation*}
where $0 \leq i, j \leq 2m -1$ and all exponents are modulo $2m$.

\begin{cor} \label{cor5}
     The generalized quaternion group $Q_m$ has a complete mapping if and only if $m$ is even.
\end{cor}

\begin{proof}
    First, we prove that if $m$ is even, then $Q_m$ has a complete mapping. Let $m = 2k$ and $|Q_{2k}| = 8k = 2^lt$, where $\gcd(2, t) = 1$ and $l \geq 3$. Suppose, for contradiction, that $H$ is a cyclic Sylow $2$-subgroup of $Q_{2k}$. By \cref{sylowsize}, $|H| = 2^l$. Since $H$ is cyclic, $H = \langle a^ib \rangle$ or $H = \langle a^i \rangle$ for some $0 \leq i \leq 4k-1$. Since $(a^ib)^2 = a^iba^ib = a^{m}$ and $a^m \neq e$, we have $|\langle a^ib \rangle| = 4 < 2^l$. Hence, $H = \langle a^i \rangle$ with $0 < i \leq 4k-1$. Let $d = \gcd(i, 4k)$. Then $\langle a^i \rangle = \langle a^d \rangle$, and since $d \mid 4k$, we have $4k = dj$ for some positive integer $j$. So, $H = \{ e, a^d, a^{2d}, \ldots, a^{(j-1)d} \}$, and hence $j = |H| = 2^l$. Thus, $8k = 2dj = 2^{l+1}d$. Since $8k = 2^lt$, we obtain $t = 2d$, which contradicts $\gcd(2, t) = 1$. Therefore, $Q_{2k}$ has no cyclic Sylow $2$-subgroup. Since $l \geq 3$, the Sylow $2$-subgroups of $Q_{2k}$ are nontrivial, and hence non-cyclic. Therefore, by \cref{thm3}, $Q_{2k}$ has a complete mapping.

    For the reverse, suppose $m$ is odd. Then $|Q_m| = 4m$, where $4$ is the largest power of $2$ dividing $4m$. Consider $\langle b \rangle = \{ e, b, a^m, a^mb \}$, which is a cyclic subgroup of order $4$. By \cref{thmlag}, every $2$-subgroup of $Q_m$ has order at most $4$, so $\langle b \rangle$ is a maximal $2$-subgroup, that is, a cyclic Sylow $2$-subgroup of $Q_m$. Therefore, by \cref{thm3}, $Q_m$ does not have a complete mapping.
\end{proof}

\section{Chromatic number of the Cayley table of commutative groups} \label{preresult}

In this section, we review known results regarding the chromatic number of the Latin square of commutative finite groups. 

By \cref{thm3}, a finite group $G$ of order $n$ has a complete mapping if and only if $L_G$ can be decomposed into $n$ disjoint transversals. Therefore, if $G$ has a complete mapping, then $\chi(L_G) = n$. 

\begin{thm} [{\cite{mahmoodian2016}}] \label{thm8}
    There exist no group $G$ of order $n$ with $\chi(L_G) = n + 1$. 
\end{thm}

\begin{proof}
    By contradiction, suppose $\chi(L_G) = n + 1$. By \cref{prop1}, $L_G$ has at least one transversal and can then be decomposed into $n$ disjoint transversals. So, $\chi(L_G) = n$ which is a contradiction with the assumption that $\chi(L_G) = n + 1$.  
\end{proof}

\begin{rem}\label{bound}
    \cref{thm3} and \cref{thm8} imply that for a group $G$ of order $n$ without complete mapping, $\chi(L_G) \geq n + 2$. 
\end{rem}

\cref{thm9} shows that the bound on $\chi(L_G)$ is tight for $G = \mathbb{Z}_{2n}$ and $G = \mathbb{Z}_{2n + 1} \times \mathbb{Z}_2$.  

\begin{thm} [{\cite{mahmoodian2016}}] \label{thm9}
    $\chi(L_{\mathbb{Z}_{2n}}) = 2n + 2$ and $\chi(L_{\mathbb{Z}_{2n + 1} \times \mathbb{Z}_2}) = 4n + 4$. 
\end{thm}

Goddyn, Halasz, and Mahmoodian \cite{Goddyn_Halasz_Mahmoodian_2019} prove \cref{thmabelian}, which completes the investigation of the chromatic number of commutative group Latin squares. 

\begin{thm}[{\cite{Goddyn_Halasz_Mahmoodian_2019}}] \label{thmabelian}
    Let $G$ be a finite commutative group of order $n$ with a cyclic Sylow 2-subgroup. Then $\chi(L_G) = n + 2$. 
\end{thm}

In \cite{Goddyn_Halasz_Mahmoodian_2019}, authors express that while there is no result regarding dihedral groups $D_n$ for odd $n$, showing that $\chi(L_{D_{n}}) = 2n + 2$ might be reachable. In the next section, we prove that $\chi(L_{D_{n}}) = 2n + 2$ for an odd number $n$.

\section{Chromatic number of $L_{D_n}$ for odd $n$}  \label{main}

While various representations of the Cayley table of $D_n$ exist, the specific structure defined below makes the coloring possible within our framework.

The dihedral group $D_n$ of order $2n$ consists of $n$ rotations $R_i = r^i$ and $n$ reflections $S_i = r^i s$ for $i \in \{0, 1, \dots, n-1\}$. 
To reveal a symmetric block structure in the Cayley table of $D_n$, we order the rows as $(R_0, \dots, R_{n-1}, S_0, \dots, S_{n-1})$ and the columns as $(R_0, \dots, R_{n-1}, S_0, S_{n-1}, \dots, S_1)$, where the reflection in the columns are indexed by $S_{-j}$ (with all indices evaluated modulo $n$).

For $i, j \in \{0, \dots, n-1\}$, the products defining the four blocks evaluate as:
\[
    \begin{aligned}
        R_i R_j &= R_{i+j}   &\quad R_i S_{-j} &= S_{i-j}   \\
        S_i R_j &= S_{i-j}   &\quad S_i S_{-j} &= R_{i+j}  
    \end{aligned}
\]
Consequently, the $2n \times 2n$ Cayley table takes the symmetric block form 
\[
    L = \begin{pmatrix} 
        A & B \\ 
        B & A 
    \end{pmatrix},
\]
where $A_{i,j} = R_{i+j}$ and $B_{i,j} = S_{i-j}$ (with indices modulo $n$). 

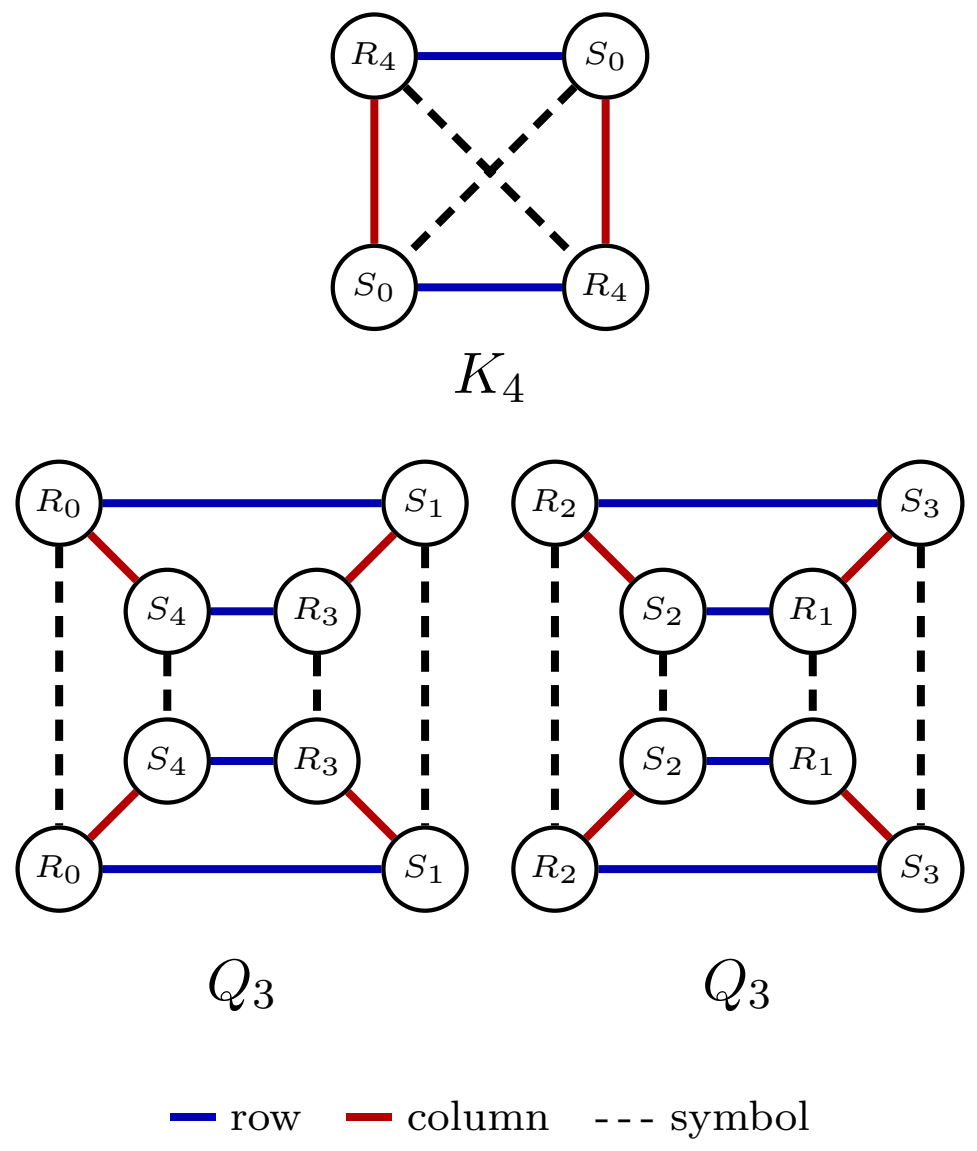
\begin{figure}[htbp]
\centering
\begin{minipage}[c]{0.54\textwidth}
\centering
\resizebox{\linewidth}{!}{%
\tiny
\setlength{\tabcolsep}{2.5pt}
\renewcommand{\arraystretch}{1.1}
\begin{tabular}{c|ccccc|ccccc}
$\cdot$ & $R_0$ & $R_1$ & $R_2$ & $R_3$ & $R_4$ & $S_0$ & $S_4$ & $S_3$ & $S_2$ & $S_1$ \\
\hline
$R_0$ & \cellcolor{green!30}$R_0$ & $R_1$ & $R_2$ & $R_3$ & $R_4$ & $S_0$ & $S_4$ & $S_3$ & $S_2$ & \cellcolor{green!30}$S_1$ \\
$R_1$ & $R_1$ & \cellcolor{green!30}$R_2$ & $R_3$ & $R_4$ & $R_0$ & $S_1$ & $S_0$ & $S_4$ & \cellcolor{green!30}$S_3$ & $S_2$ \\
$R_2$ & $R_2$ & $R_3$ & \cellcolor{green!100}$\mathbf{R_4}$ & $R_0$ & $R_1$ & $S_2$ & $S_1$ & \cellcolor{green!100}$\mathbf{S_0}$ & $S_4$ & $S_3$ \\
$R_3$ & $R_3$ & $R_4$ & $R_0$ & \cellcolor{green!30}$R_1$ & $R_2$ & $S_3$ & \cellcolor{green!30}$S_2$ & $S_1$ & $S_0$ & $S_4$ \\
$R_4$ & $R_4$ & $R_0$ & $R_1$ & $R_2$ & \cellcolor{green!30}$R_3$ & \cellcolor{green!30}$S_4$ & $S_3$ & $S_2$ & $S_1$ & $S_0$ \\
\hline
$S_0$ & $S_0$ & $S_4$ & $S_3$ & $S_2$ & \cellcolor{green!30}$S_1$ & \cellcolor{green!30}$R_0$ & $R_1$ & $R_2$ & $R_3$ & $R_4$ \\
$S_1$ & $S_1$ & $S_0$ & $S_4$ & \cellcolor{green!30}$S_3$ & $S_2$ & $R_1$ & \cellcolor{green!30}$R_2$ & $R_3$ & $R_4$ & $R_0$ \\
$S_2$ & $S_2$ & $S_1$ & \cellcolor{green!100}$\mathbf{S_0}$ & $S_4$ & $S_3$ & $R_2$ & $R_3$ & \cellcolor{green!100}$\mathbf{R_4}$ & $R_0$ & $R_1$ \\
$S_3$ & $S_3$ & \cellcolor{green!30}$S_2$ & $S_1$ & $S_0$ & $S_4$ & $R_3$ & $R_4$ & $R_0$ & \cellcolor{green!30}$R_1$ & $R_2$ \\
$S_4$ & \cellcolor{green!30}$S_4$ & $S_3$ & $S_2$ & $S_1$ & $S_0$ & $R_4$ & $R_0$ & $R_1$ & $R_2$ & \cellcolor{green!30}$R_3$ \\
\end{tabular}}
\end{minipage}\hfill
\begin{minipage}[c]{0.43\textwidth}
\centering
\begin{tikzpicture}[scale=0.85]
  \node[vtx] (a) at (0,1.8) {$R_4$};
  \node[vtx] (b) at (1.8,1.8) {$S_0$};
  \node[vtx] (c) at (1.8,0) {$R_4$};
  \node[vtx] (d) at (0,0) {$S_0$};
  \begin{scope}[on background layer]
    \draw[redge] (a)--(b) (d)--(c);
    \draw[cedge] (a)--(d) (b)--(c);
    \draw[sedge] (a)--(c) (b)--(d);
  \end{scope}
  \node[font=\small] at (0.9,-0.7) {$K_4$};
\end{tikzpicture}

\vspace{6pt}
\qcube{$R_0$}{$S_1$}{$S_1$}{$R_0$}{$S_4$}{$R_3$}{$R_3$}{$S_4$}\hspace{8pt}%
\qcube{$R_2$}{$S_3$}{$S_3$}{$R_2$}{$S_2$}{$R_1$}{$R_1$}{$S_2$}

\vspace{4pt}
{\scriptsize
\textcolor{blue!70!black}{\rule[0.5ex]{1em}{1.5pt}} row\quad
\textcolor{red!70!black}{\rule[0.5ex]{1em}{1.5pt}} column\quad
-\,-\,- symbol}
\end{minipage}
\caption{The $10 \times 10$ Cayley table for $D_5$ with the 2-plex $\mathcal{P}_{(0,0)}$ highlighted (left), and the graph it induces (right), with vertices labeled by their symbols.} 
\label{fig:symmetric_d5_2plex}
\end{figure}

Throughout this section, $L$ denotes the $2n \times 2n$ Cayley table with the block structure defined above. See Figures~\ref{fig:symmetric_d5_2plex} and~\ref{fig:d9_UV} for the Cayley tables of $D_5$ and $D_9$, respectively.

\begin{lem}\label{lem:graph-dec}
Let $n \ge 3$ be an odd integer, and let $L$ be the  Cayley table of $D_n$.  For $a,b\in \mathbb{Z}_n$, let
\[
\begin{aligned}
\mathcal{P}_{(a,b)}
={}& \{(i+a, \, i; \, R_{2i+a}) \mid i \in \mathbb{Z}_n\} \\
&\cup \{(i+b, \, 2n-1-i; \, S_{2i+1+b}) \mid i \in \mathbb{Z}_n\} \\
&\cup \{(n+i+b, \, n-1-i; \, S_{2i+1+b}) \mid i \in \mathbb{Z}_n\} \\
&\cup \{(n+i+a, \, n+i; \, R_{2i+a}) \mid i \in \mathbb{Z}_n\}.
\end{aligned}
\]
Then $\mathcal{P}_{(a,b)}$ is a $2$-plex with induced graph 
\[
G_{(a,b)} \;\cong\; K_4 \;\oplus\; \left(\frac{n - 1}{2}\right) Q_3.
\]
\end{lem}

\begin{proof}
We partition the cells of $\mathcal{P}_{(a,b)}$ into four families; for $i \in \mathbb{Z}_n$, let
\begin{align*}
    u_i &:= (i + a, \; i; \; R_{2i + a})\in A_{\text{top-left}}, &
    x_i &:= (i + b, \; 2n - 1 - i; \; S_{2i + 1 + b})\in B_{\text{top-right}}, \\
    y_i &:= (n + i + b, \; n - 1 - i; \; S_{2i + 1 + b})\in B_{\text{bottom-left}}, &
    v_i &:= (n + i + a, \; n + i; \; R_{2i + a})\in A_{\text{bottom-right}}. 
\end{align*}
Note that the index arithmetic is evaluated modulo $n$. 

Since $\gcd(2, n) = 1$, the maps $i \mapsto 2i + a \pmod n$ and $i \mapsto 2i + 1 + b \pmod n$ are bijections on $\mathbb{Z}_n$. Consequently, every row, column, and group symbol appears in precisely two cells, confirming that $\mathcal{P}_{(a,b)}$ is a $2$-plex and its induced graph $G_{(a,b)}$ is $3$-regular. Thus, the edge set of $G_{(a,b)}$ decomposes into the row, column, and symbol matchings: $E = M_{\text{row}} \dot\cup M_{\text{col}} \dot\cup M_{\text{sym}}$. 

Explicitly, for all $i \in \mathbb{Z}_n$, these matchings pair the vertices as:
\begin{align*}
    M_{\text{row}} &: \quad u_i \sim x_{i + a - b} \quad \text{and} \quad y_i \sim v_{i + b - a}, \\
    M_{\text{col}} &: \quad u_i \sim y_{n - 1 - i} \quad \text{and} \quad x_i \sim v_{n - 1 - i}, \\
    M_{\text{sym}} &: \quad u_i \sim v_i \quad \text{and} \quad x_i \sim y_i.
\end{align*}
The alternating union $M_{\text{row}} \cup M_{\text{col}}$ partitions the vertex set into $n$ disjoint $4$-cycles given by
\[
C_4(i) = \bigl(u_i, \, x_{i + a - b}, \, v_{n - 1 - i - a + b}, \, y_{n - 1 - i}\bigr).
\]
Under $M_{\text{sym}}$, the vertex $u_i$ is adjacent to $v_i$. For any $j \in \mathbb{Z}_n$,
\[
C_4(j) = \bigl(u_j, \, x_{j + a - b}, \, v_{n - 1 - j - a + b}, \, y_{n - 1 - j}\bigr),
\]
so $v_i$ is the third vertex of $C_4(j)$ precisely when $n - 1 - j - a + b \equiv i$, that is, $j = b - a - 1 - i$. Thus, let $\sigma(i) = b - a - 1 - i$; since $\sigma(\sigma(i)) = i$, the map $\sigma$ is an involution on $\mathbb{Z}_n$. Substituting $j = \sigma(i)$ gives
\[
C_4(\sigma(i)) = \bigl(u_{\sigma(i)}, \, x_{n-1-i}, \, v_i, \, y_{i+a-b}\bigr).
\]
Comparing $C_4(i)$ and $C_4(\sigma(i))$, the symbol edges $u_i \sim v_i$, $x_{i+a-b} \sim y_{i+a-b}$, $v_{\sigma(i)} \sim u_{\sigma(i)}$, and $y_{n-1-i} \sim x_{n-1-i}$ join the first, second, third, and fourth vertices of $C_4(i)$ to the third, fourth, first, and second vertices of $C_4(\sigma(i))$, respectively. Because $n$ is odd, $\sigma$ has a unique fixed point at $i^* = 2^{-1}(b - a - 1)$, where $2^{-1}$ denotes the inverse of $2$ modulo $n$. At this fixed point, the symbol chords join opposite vertices of $C_4(i^*)$ and convert it into a complete graph $K_4$. The remaining $n - 1$ indices partition into $(n - 1)/2$ disjoint pairs $\{i, \sigma(i)\}$. Each index pair defines a pair of distinct $4$-cycles, $C_4(i)$ and $C_4(\sigma(i))$, connected under $M_{\text{sym}}$ by a cyclic shift of two positions; since this shift preserves adjacency, they form an isomorphic copy of the $3$-dimensional hypercube $Q_3 \cong C_4 \times K_2$. Thus, $G_{(a,b)}$ decomposes into one $K_4$ and $(n - 1)/2$ copies of $Q_3$. For an example, see \cref{fig:symmetric_d5_2plex}.
\end{proof}


\begin{lem} \label{lem2-plex}
Let $n \ge 5$ be an odd integer and let $L$ be the Cayley table of $D_n$. There exist pairs $(a_k, b_k) \in \mathbb{Z}_n^2$, $0 \le k \le n-1$, such that the $2$-plexes $\mathcal{P}_{(a_k, b_k)}$ (defined as in~\cref{lem:graph-dec}) jointly partition $L$, and every row and column of $L$ meets the $K_4$ component of at most two of these $2$-plexes.
\end{lem}

\begin{proof}
By the proof of \cref{lem:graph-dec}, the $K_4$ component of $\mathcal{P}_{(a_k,b_k)}$ is a $2\times 2$ rectangle of cells. For $i,j\in\mathbb{Z}_n$, its four cells are
\begin{align*}
    &((i + a_k) \bmod n, \, i), && ((j + b_k) \bmod n, \, 2n - 1 - j), \\
    &(n + (j + b_k) \bmod n, \, n - 1 - j), && (n + (i + a_k) \bmod n, \, n + i),
\end{align*}
which are the top-left, top-right, bottom-left and bottom-right cells, respectively. The top cells share a row when $j \equiv i + a_k - b_k \pmod n$, and the left cells share a column when $j = n-1-i$. Thus, $2i \equiv (n-1) + (b_k - a_k) \pmod n$.
Since $n$ is odd, $2$ is invertible modulo $n$ and thus we can determines the unique column coordinate, $\tfrac{n-1}{2} + 2^{-1}(b_k - a_k),$ and the unique row coordinate 
$\tfrac{n-1}{2} + 2^{-1}(a_k + b_k)$.

Thus, the maximum number of $K_4$ components intersecting any row or column of $L$ is bounded by the maximum frequencies of differences and sums in $\mathbb{Z}_n$:
\begin{equation} \label{max2}
D = \max_{d \in \mathbb{Z}_n} \big| \{k \in \mathbb{Z}_n \mid b_k - a_k \equiv d \pmod n\} \big|, \quad 
S = \max_{s \in \mathbb{Z}_n} \big| \{k \in \mathbb{Z}_n \mid a_k + b_k \equiv s \pmod n\} \big|.
\end{equation}

Let $a_k = k$ for $0 \le k \le n-1$, so $(a_k)$ is a permutation of $\mathbb{Z}_n$. It remains to choose $(b_k)$ such that two conditions are satisfied. First, $(b_k)$ must itself be a permutation of $\mathbb{Z}_n$, ensuring that the $n$ 2-plexes $\mathcal{P}_{(a_k, b_k)}$ partition $L$. Second, $(b_k)$ must satisfy $D, S \le 2$, where $D$ and $S$ are defined  in~\eqref{max2}. We consider two cases.

\medskip
\noindent\textbf{Case 1: $3 \nmid n$.}
Let $b_k = 2k \bmod n$. Since $\gcd(2,n)=1$, $(b_k)$ is a permutation of $\mathbb{Z}_n$. Moreover, $b_k - a_k \equiv k$ and $a_k + b_k \equiv 3k \pmod n$. Since $\gcd(3,n)=1$, both $k \mapsto k$ and $k \mapsto 3k$ are bijections of $\mathbb{Z}_n$, so $D = S = 1$. Hence every row and column of $L$ meets the $K_4$ component of exactly one of the $2$-plexes.

\medskip
\noindent\textbf{Case 2: $3 \mid n$.} 
Write $n = 3m$ with $m \ge 3$ odd. Partition $\mathbb{Z}_n$ into three intervals: $B_1 = \{0, \dots, m-1\}$, $B_2 = \{m, \dots, 2m-1\}$, and $B_3 = \{2m, \dots, 3m-1\}$, and let $B = B_2 \cup B_3$. 

If we choose $b_k = 2k \pmod n$ for all $k$, then $a_k + b_k \equiv 3k \pmod n$ takes every multiple of $3$ three times, namely at $k$, $k+m$ and $k+2m$ (one index from each of $B_1$, $B_2$, $B_3$), so $S = 3$. To avoid this, we keep $b_k = 2k \pmod n$ on $B$ and modify it only on $B_1$. Let $t = 2$ if $m \equiv 1 \pmod 3$ and $t = 1$ otherwise, so that $t \not\equiv 0,\, m \pmod 3$, and define
\[
b_k = \begin{cases}
2\big((k+t) \bmod m\big), & k \in B_1, \\
2k \pmod n, & k \in B.
\end{cases}
\]
Explicitly, on $B_1$ this evaluates as
\[
b_k = \begin{cases}
2k + 2t, & k < m-t, \\
2k + 2t - 2m, & m-t \le k < m.
\end{cases}
\]

We first verify that $(b_k)$ is a permutation of $\mathbb{Z}_n$. On $B_1$, $k \mapsto b_k$ is a bijection onto $\{0, 2, \dots, 2m-2\}$, which is the image of $B_1$ under the bijection $k \mapsto 2k \pmod n$ of $\mathbb{Z}_n$. Hence $b_k = 2k \pmod n$ on $B$ takes exactly the remaining values, and $(b_k)$ is a permutation of $\mathbb{Z}_n$.

Next, we show that $b_k - a_k$ takes each value in $\mathbb{Z}_n$ at most twice, so $D \le 2$.
First, note that $b_k - a_k = k$ is injective on $B$. On $B_1$, using $-2m \equiv m \pmod n$, we have
\[
b_k - a_k \equiv \begin{cases}
k + 2t \pmod n \ \in [\,2t,\; m+t-1\,], & k < m-t, \\
k + 2t + m \pmod n \ \in [\,2m+t,\; 2m+2t-1\,], & m-t \le k < m.
\end{cases}
\]

These two intervals are disjoint subsets of $\{0, \dots, n-1\}$, meaning $b_k - a_k$ is also injective on $B_1$. Consequently, every residue occurs at most once from $B$ and at most once from $B_1$, giving $D \le 2$.

Finally, we show that $a_k + b_k$ takes each value in $\mathbb{Z}_n$ at most twice, so $S \le 2$. On $B$, the sum $a_k + b_k \equiv 3k \pmod n$ covers each multiple of $3$ exactly twice and no other values. On $B_1$, the sums are given by
\[
a_k + b_k \equiv \begin{cases}
3k + 2t \pmod n, & k < m-t, \\
3k + 2(t-m) \pmod n, & m-t \le k < m.
\end{cases}
\]
By our choice of $t$, neither $2t$ nor $2(t-m)$ is divisible by $3$. Consequently, none of these sums is a multiple of $3$, so sums from $B_1$ never coincide with sums from $B$. Furthermore, within each sub-interval,  each sum occurs at most once per sub-interval, and thus at most twice across all of $B_1$. 

Therefore, $S \le 2$, which completes the proof.
\Cref{tab:d9_2plexes} shows Case~2 for $n=9$, where $m=3$ and $t=1$. Here $(a_k)_{k=0}^{8}=k$ and $(b_k)_{k=0}^{8}=(2,4,0,6,8,1,3,5,7)$. Each $2$-plex is drawn in its own light color, and its $K_4$ component is shown in a darker shade of the same color. As the figure shows, every row and column of $L$ meets at most two of the $K_4$ components.
\end{proof}

\begin{table}[h]
\centering
\tiny
\setlength{\tabcolsep}{2.5pt}
\renewcommand{\arraystretch}{1.1}
\resizebox{\textwidth}{!}{%
\begin{tabular}{c|ccccccccc|ccccccccc}
$\cdot$ & $R_{0}$ & $R_{1}$ & $R_{2}$ & $R_{3}$ & $R_{4}$ & $R_{5}$ & $R_{6}$ & $R_{7}$ & $R_{8}$ & $S_{0}$ & $S_{8}$ & $S_{7}$ & $S_{6}$ & $S_{5}$ & $S_{4}$ & $S_{3}$ & $S_{2}$ & $S_{1}$ \\
\hline
$R_{0}$ & \cellcolor{pl0!28}$R_{0}$ & \cellcolor{pl8!28}$R_{1}$ & \cellcolor{pl7!28}$R_{2}$ & \cellcolor{pl6!28}$R_{3}$ & \cellcolor{pl5!28}$R_{4}$ & \cellcolor{pl4!28}$R_{5}$ & \cellcolor{pl3!28}$R_{6}$ & \cellcolor{pl2!28}$R_{7}$ & \cellcolor{pl1!28}$R_{8}$ & \cellcolor{pl5!28}$S_{0}$ & \cellcolor{pl0!28}$S_{8}$ & \cellcolor{pl6!28}$S_{7}$ & \cellcolor{pl1!28}$S_{6}$ & \cellcolor{pl7!28}$S_{5}$ & \cellcolor{pl3!28}$S_{4}$ & \cellcolor{pl8!28}$S_{3}$ & \cellcolor{pl4!28}$S_{2}$ & \cellcolor{pl2!28}$S_{1}$ \\
$R_{1}$ & \cellcolor{pl1!28}$R_{1}$ & \cellcolor{pl0!28}$R_{2}$ & \cellcolor{pl8!28}$R_{3}$ & \cellcolor{pl7!75}$\mathbf{R_{4}}$ & \cellcolor{pl6!28}$R_{5}$ & \cellcolor{pl5!28}$R_{6}$ & \cellcolor{pl4!75}$\mathbf{R_{7}}$ & \cellcolor{pl3!28}$R_{8}$ & \cellcolor{pl2!28}$R_{0}$ & \cellcolor{pl0!28}$S_{1}$ & \cellcolor{pl6!28}$S_{0}$ & \cellcolor{pl1!28}$S_{8}$ & \cellcolor{pl7!75}$\mathbf{S_{7}}$ & \cellcolor{pl3!28}$S_{6}$ & \cellcolor{pl8!28}$S_{5}$ & \cellcolor{pl4!75}$\mathbf{S_{4}}$ & \cellcolor{pl2!28}$S_{3}$ & \cellcolor{pl5!28}$S_{2}$ \\
$R_{2}$ & \cellcolor{pl2!28}$R_{2}$ & \cellcolor{pl1!75}$\mathbf{R_{3}}$ & \cellcolor{pl0!28}$R_{4}$ & \cellcolor{pl8!28}$R_{5}$ & \cellcolor{pl7!28}$R_{6}$ & \cellcolor{pl6!28}$R_{7}$ & \cellcolor{pl5!28}$R_{8}$ & \cellcolor{pl4!28}$R_{0}$ & \cellcolor{pl3!28}$R_{1}$ & \cellcolor{pl6!28}$S_{2}$ & \cellcolor{pl1!75}$\mathbf{S_{1}}$ & \cellcolor{pl7!28}$S_{0}$ & \cellcolor{pl3!28}$S_{8}$ & \cellcolor{pl8!28}$S_{7}$ & \cellcolor{pl4!28}$S_{6}$ & \cellcolor{pl2!28}$S_{5}$ & \cellcolor{pl5!28}$S_{4}$ & \cellcolor{pl0!28}$S_{3}$ \\
$R_{3}$ & \cellcolor{pl3!28}$R_{3}$ & \cellcolor{pl2!28}$R_{4}$ & \cellcolor{pl1!28}$R_{5}$ & \cellcolor{pl0!28}$R_{6}$ & \cellcolor{pl8!28}$R_{7}$ & \cellcolor{pl7!28}$R_{8}$ & \cellcolor{pl6!28}$R_{0}$ & \cellcolor{pl5!28}$R_{1}$ & \cellcolor{pl4!28}$R_{2}$ & \cellcolor{pl1!28}$S_{3}$ & \cellcolor{pl7!28}$S_{2}$ & \cellcolor{pl3!28}$S_{1}$ & \cellcolor{pl8!28}$S_{0}$ & \cellcolor{pl4!28}$S_{8}$ & \cellcolor{pl2!28}$S_{7}$ & \cellcolor{pl5!28}$S_{6}$ & \cellcolor{pl0!28}$S_{5}$ & \cellcolor{pl6!28}$S_{4}$ \\
$R_{4}$ & \cellcolor{pl4!28}$R_{4}$ & \cellcolor{pl3!75}$\mathbf{R_{5}}$ & \cellcolor{pl2!28}$R_{6}$ & \cellcolor{pl1!28}$R_{7}$ & \cellcolor{pl0!28}$R_{8}$ & \cellcolor{pl8!28}$R_{0}$ & \cellcolor{pl7!28}$R_{1}$ & \cellcolor{pl6!75}$\mathbf{R_{2}}$ & \cellcolor{pl5!28}$R_{3}$ & \cellcolor{pl7!28}$S_{4}$ & \cellcolor{pl3!75}$\mathbf{S_{3}}$ & \cellcolor{pl8!28}$S_{2}$ & \cellcolor{pl4!28}$S_{1}$ & \cellcolor{pl2!28}$S_{0}$ & \cellcolor{pl5!28}$S_{8}$ & \cellcolor{pl0!28}$S_{7}$ & \cellcolor{pl6!75}$\mathbf{S_{6}}$ & \cellcolor{pl1!28}$S_{5}$ \\
$R_{5}$ & \cellcolor{pl5!28}$R_{5}$ & \cellcolor{pl4!28}$R_{6}$ & \cellcolor{pl3!28}$R_{7}$ & \cellcolor{pl2!75}$\mathbf{R_{8}}$ & \cellcolor{pl1!28}$R_{0}$ & \cellcolor{pl0!75}$\mathbf{R_{1}}$ & \cellcolor{pl8!28}$R_{2}$ & \cellcolor{pl7!28}$R_{3}$ & \cellcolor{pl6!28}$R_{4}$ & \cellcolor{pl3!28}$S_{5}$ & \cellcolor{pl8!28}$S_{4}$ & \cellcolor{pl4!28}$S_{3}$ & \cellcolor{pl2!75}$\mathbf{S_{2}}$ & \cellcolor{pl5!28}$S_{1}$ & \cellcolor{pl0!75}$\mathbf{S_{0}}$ & \cellcolor{pl6!28}$S_{8}$ & \cellcolor{pl1!28}$S_{7}$ & \cellcolor{pl7!28}$S_{6}$ \\
$R_{6}$ & \cellcolor{pl6!28}$R_{6}$ & \cellcolor{pl5!28}$R_{7}$ & \cellcolor{pl4!28}$R_{8}$ & \cellcolor{pl3!28}$R_{0}$ & \cellcolor{pl2!28}$R_{1}$ & \cellcolor{pl1!28}$R_{2}$ & \cellcolor{pl0!28}$R_{3}$ & \cellcolor{pl8!28}$R_{4}$ & \cellcolor{pl7!28}$R_{5}$ & \cellcolor{pl8!28}$S_{6}$ & \cellcolor{pl4!28}$S_{5}$ & \cellcolor{pl2!28}$S_{4}$ & \cellcolor{pl5!28}$S_{3}$ & \cellcolor{pl0!28}$S_{2}$ & \cellcolor{pl6!28}$S_{1}$ & \cellcolor{pl1!28}$S_{0}$ & \cellcolor{pl7!28}$S_{8}$ & \cellcolor{pl3!28}$S_{7}$ \\
$R_{7}$ & \cellcolor{pl7!28}$R_{7}$ & \cellcolor{pl6!28}$R_{8}$ & \cellcolor{pl5!75}$\mathbf{R_{0}}$ & \cellcolor{pl4!28}$R_{1}$ & \cellcolor{pl3!28}$R_{2}$ & \cellcolor{pl2!28}$R_{3}$ & \cellcolor{pl1!28}$R_{4}$ & \cellcolor{pl0!28}$R_{5}$ & \cellcolor{pl8!75}$\mathbf{R_{6}}$ & \cellcolor{pl4!28}$S_{7}$ & \cellcolor{pl2!28}$S_{6}$ & \cellcolor{pl5!75}$\mathbf{S_{5}}$ & \cellcolor{pl0!28}$S_{4}$ & \cellcolor{pl6!28}$S_{3}$ & \cellcolor{pl1!28}$S_{2}$ & \cellcolor{pl7!28}$S_{1}$ & \cellcolor{pl3!28}$S_{0}$ & \cellcolor{pl8!75}$\mathbf{S_{8}}$ \\
$R_{8}$ & \cellcolor{pl8!28}$R_{8}$ & \cellcolor{pl7!28}$R_{0}$ & \cellcolor{pl6!28}$R_{1}$ & \cellcolor{pl5!28}$R_{2}$ & \cellcolor{pl4!28}$R_{3}$ & \cellcolor{pl3!28}$R_{4}$ & \cellcolor{pl2!28}$R_{5}$ & \cellcolor{pl1!28}$R_{6}$ & \cellcolor{pl0!28}$R_{7}$ & \cellcolor{pl2!28}$S_{8}$ & \cellcolor{pl5!28}$S_{7}$ & \cellcolor{pl0!28}$S_{6}$ & \cellcolor{pl6!28}$S_{5}$ & \cellcolor{pl1!28}$S_{4}$ & \cellcolor{pl7!28}$S_{3}$ & \cellcolor{pl3!28}$S_{2}$ & \cellcolor{pl8!28}$S_{1}$ & \cellcolor{pl4!28}$S_{0}$ \\
\hline
$S_{0}$ & \cellcolor{pl5!28}$S_{0}$ & \cellcolor{pl0!28}$S_{8}$ & \cellcolor{pl6!28}$S_{7}$ & \cellcolor{pl1!28}$S_{6}$ & \cellcolor{pl7!28}$S_{5}$ & \cellcolor{pl3!28}$S_{4}$ & \cellcolor{pl8!28}$S_{3}$ & \cellcolor{pl4!28}$S_{2}$ & \cellcolor{pl2!28}$S_{1}$ & \cellcolor{pl0!28}$R_{0}$ & \cellcolor{pl8!28}$R_{1}$ & \cellcolor{pl7!28}$R_{2}$ & \cellcolor{pl6!28}$R_{3}$ & \cellcolor{pl5!28}$R_{4}$ & \cellcolor{pl4!28}$R_{5}$ & \cellcolor{pl3!28}$R_{6}$ & \cellcolor{pl2!28}$R_{7}$ & \cellcolor{pl1!28}$R_{8}$ \\
$S_{1}$ & \cellcolor{pl0!28}$S_{1}$ & \cellcolor{pl6!28}$S_{0}$ & \cellcolor{pl1!28}$S_{8}$ & \cellcolor{pl7!75}$\mathbf{S_{7}}$ & \cellcolor{pl3!28}$S_{6}$ & \cellcolor{pl8!28}$S_{5}$ & \cellcolor{pl4!75}$\mathbf{S_{4}}$ & \cellcolor{pl2!28}$S_{3}$ & \cellcolor{pl5!28}$S_{2}$ & \cellcolor{pl1!28}$R_{1}$ & \cellcolor{pl0!28}$R_{2}$ & \cellcolor{pl8!28}$R_{3}$ & \cellcolor{pl7!75}$\mathbf{R_{4}}$ & \cellcolor{pl6!28}$R_{5}$ & \cellcolor{pl5!28}$R_{6}$ & \cellcolor{pl4!75}$\mathbf{R_{7}}$ & \cellcolor{pl3!28}$R_{8}$ & \cellcolor{pl2!28}$R_{0}$ \\
$S_{2}$ & \cellcolor{pl6!28}$S_{2}$ & \cellcolor{pl1!75}$\mathbf{S_{1}}$ & \cellcolor{pl7!28}$S_{0}$ & \cellcolor{pl3!28}$S_{8}$ & \cellcolor{pl8!28}$S_{7}$ & \cellcolor{pl4!28}$S_{6}$ & \cellcolor{pl2!28}$S_{5}$ & \cellcolor{pl5!28}$S_{4}$ & \cellcolor{pl0!28}$S_{3}$ & \cellcolor{pl2!28}$R_{2}$ & \cellcolor{pl1!75}$\mathbf{R_{3}}$ & \cellcolor{pl0!28}$R_{4}$ & \cellcolor{pl8!28}$R_{5}$ & \cellcolor{pl7!28}$R_{6}$ & \cellcolor{pl6!28}$R_{7}$ & \cellcolor{pl5!28}$R_{8}$ & \cellcolor{pl4!28}$R_{0}$ & \cellcolor{pl3!28}$R_{1}$ \\
$S_{3}$ & \cellcolor{pl1!28}$S_{3}$ & \cellcolor{pl7!28}$S_{2}$ & \cellcolor{pl3!28}$S_{1}$ & \cellcolor{pl8!28}$S_{0}$ & \cellcolor{pl4!28}$S_{8}$ & \cellcolor{pl2!28}$S_{7}$ & \cellcolor{pl5!28}$S_{6}$ & \cellcolor{pl0!28}$S_{5}$ & \cellcolor{pl6!28}$S_{4}$ & \cellcolor{pl3!28}$R_{3}$ & \cellcolor{pl2!28}$R_{4}$ & \cellcolor{pl1!28}$R_{5}$ & \cellcolor{pl0!28}$R_{6}$ & \cellcolor{pl8!28}$R_{7}$ & \cellcolor{pl7!28}$R_{8}$ & \cellcolor{pl6!28}$R_{0}$ & \cellcolor{pl5!28}$R_{1}$ & \cellcolor{pl4!28}$R_{2}$ \\
$S_{4}$ & \cellcolor{pl7!28}$S_{4}$ & \cellcolor{pl3!75}$\mathbf{S_{3}}$ & \cellcolor{pl8!28}$S_{2}$ & \cellcolor{pl4!28}$S_{1}$ & \cellcolor{pl2!28}$S_{0}$ & \cellcolor{pl5!28}$S_{8}$ & \cellcolor{pl0!28}$S_{7}$ & \cellcolor{pl6!75}$\mathbf{S_{6}}$ & \cellcolor{pl1!28}$S_{5}$ & \cellcolor{pl4!28}$R_{4}$ & \cellcolor{pl3!75}$\mathbf{R_{5}}$ & \cellcolor{pl2!28}$R_{6}$ & \cellcolor{pl1!28}$R_{7}$ & \cellcolor{pl0!28}$R_{8}$ & \cellcolor{pl8!28}$R_{0}$ & \cellcolor{pl7!28}$R_{1}$ & \cellcolor{pl6!75}$\mathbf{R_{2}}$ & \cellcolor{pl5!28}$R_{3}$ \\
$S_{5}$ & \cellcolor{pl3!28}$S_{5}$ & \cellcolor{pl8!28}$S_{4}$ & \cellcolor{pl4!28}$S_{3}$ & \cellcolor{pl2!75}$\mathbf{S_{2}}$ & \cellcolor{pl5!28}$S_{1}$ & \cellcolor{pl0!75}$\mathbf{S_{0}}$ & \cellcolor{pl6!28}$S_{8}$ & \cellcolor{pl1!28}$S_{7}$ & \cellcolor{pl7!28}$S_{6}$ & \cellcolor{pl5!28}$R_{5}$ & \cellcolor{pl4!28}$R_{6}$ & \cellcolor{pl3!28}$R_{7}$ & \cellcolor{pl2!75}$\mathbf{R_{8}}$ & \cellcolor{pl1!28}$R_{0}$ & \cellcolor{pl0!75}$\mathbf{R_{1}}$ & \cellcolor{pl8!28}$R_{2}$ & \cellcolor{pl7!28}$R_{3}$ & \cellcolor{pl6!28}$R_{4}$ \\
$S_{6}$ & \cellcolor{pl8!28}$S_{6}$ & \cellcolor{pl4!28}$S_{5}$ & \cellcolor{pl2!28}$S_{4}$ & \cellcolor{pl5!28}$S_{3}$ & \cellcolor{pl0!28}$S_{2}$ & \cellcolor{pl6!28}$S_{1}$ & \cellcolor{pl1!28}$S_{0}$ & \cellcolor{pl7!28}$S_{8}$ & \cellcolor{pl3!28}$S_{7}$ & \cellcolor{pl6!28}$R_{6}$ & \cellcolor{pl5!28}$R_{7}$ & \cellcolor{pl4!28}$R_{8}$ & \cellcolor{pl3!28}$R_{0}$ & \cellcolor{pl2!28}$R_{1}$ & \cellcolor{pl1!28}$R_{2}$ & \cellcolor{pl0!28}$R_{3}$ & \cellcolor{pl8!28}$R_{4}$ & \cellcolor{pl7!28}$R_{5}$ \\
$S_{7}$ & \cellcolor{pl4!28}$S_{7}$ & \cellcolor{pl2!28}$S_{6}$ & \cellcolor{pl5!75}$\mathbf{S_{5}}$ & \cellcolor{pl0!28}$S_{4}$ & \cellcolor{pl6!28}$S_{3}$ & \cellcolor{pl1!28}$S_{2}$ & \cellcolor{pl7!28}$S_{1}$ & \cellcolor{pl3!28}$S_{0}$ & \cellcolor{pl8!75}$\mathbf{S_{8}}$ & \cellcolor{pl7!28}$R_{7}$ & \cellcolor{pl6!28}$R_{8}$ & \cellcolor{pl5!75}$\mathbf{R_{0}}$ & \cellcolor{pl4!28}$R_{1}$ & \cellcolor{pl3!28}$R_{2}$ & \cellcolor{pl2!28}$R_{3}$ & \cellcolor{pl1!28}$R_{4}$ & \cellcolor{pl0!28}$R_{5}$ & \cellcolor{pl8!75}$\mathbf{R_{6}}$ \\
$S_{8}$ & \cellcolor{pl2!28}$S_{8}$ & \cellcolor{pl5!28}$S_{7}$ & \cellcolor{pl0!28}$S_{6}$ & \cellcolor{pl6!28}$S_{5}$ & \cellcolor{pl1!28}$S_{4}$ & \cellcolor{pl7!28}$S_{3}$ & \cellcolor{pl3!28}$S_{2}$ & \cellcolor{pl8!28}$S_{1}$ & \cellcolor{pl4!28}$S_{0}$ & \cellcolor{pl8!28}$R_{8}$ & \cellcolor{pl7!28}$R_{0}$ & \cellcolor{pl6!28}$R_{1}$ & \cellcolor{pl5!28}$R_{2}$ & \cellcolor{pl4!28}$R_{3}$ & \cellcolor{pl3!28}$R_{4}$ & \cellcolor{pl2!28}$R_{5}$ & \cellcolor{pl1!28}$R_{6}$ & \cellcolor{pl0!28}$R_{7}$ \\
\end{tabular}}

\vspace{15pt}
{\footnotesize
\setlength{\tabcolsep}{4pt}
\renewcommand{\arraystretch}{1.3}
\begin{tabular}{!{\vrule width 1.2pt}l!{\vrule width 1.2pt}*{9}{>{\centering\arraybackslash}p{1.6em}}!{\vrule width 1.2pt}}
\noalign{\hrule height 1.2pt}
$2$-plex & \cellcolor{pl0!28} & \cellcolor{pl1!28} & \cellcolor{pl2!28} & \cellcolor{pl3!28} & \cellcolor{pl4!28} & \cellcolor{pl5!28} & \cellcolor{pl6!28} & \cellcolor{pl7!28} & \cellcolor{pl8!28} \\
$K_4$ & \cellcolor{pl0!75} & \cellcolor{pl1!75} & \cellcolor{pl2!75} & \cellcolor{pl3!75} & \cellcolor{pl4!75} & \cellcolor{pl5!75} & \cellcolor{pl6!75} & \cellcolor{pl7!75} & \cellcolor{pl8!75} \\
\noalign{\hrule height 1.2pt}
$a_k$ & $0$ & $1$ & $2$ & $3$ & $4$ & $5$ & $6$ & $7$ & $8$ \\
$b_k$ & $2$ & $4$ & $0$ & $6$ & $8$ & $1$ & $3$ & $5$ & $7$ \\
\noalign{\hrule height 1.2pt}
\end{tabular}}
\caption{The nine $2$-plexes $\mathcal{P}_{(a_k,b_k)}$ of $D_9$, colored as in the key below the table; the darker cells form the $K_4$ component of each $2$-plex.}
\label{tab:d9_2plexes}
\end{table}


We are now ready to present and prove our main theorem using \cref{lem:graph-dec} and \cref{lem2-plex}.

\begin{thm} \label{thmmain}
Let $n \ge 3$ be odd. Then $\chi(L_{D_n}) = 2n+2$.
\end{thm}

\begin{proof}
By \cref{bound}, $\chi(L_{D_n}) \geq 2n+2$, so it suffices to give a proper coloring of $L_{D_n}$ with $2n+2$ colors.

For $n = 3$, the Cayley table $L_{D_3}$ has order $6$. By \cref{bound}, $\chi(L_{D_3}) \geq 8$. On the other hand, since \cref{eq1} holds for every Latin square of order at most $8$~\cite{mahmoodian2016}, we have $\chi(L_{D_3}) \leq 6 + 2 = 8$. Thus $\chi(L_{D_3}) = 8$. Hence we may assume $n \ge 5$.

Take the pairs $(a_k,b_k)$ from \cref{lem2-plex} and write $\mathcal{P}_k = \mathcal{P}_{(a_k,b_k)}$. By \cref{lem2-plex}, the $\mathcal{P}_k$ partition the cells of $L_{D_n}$, and by \cref{lem:graph-dec} the induced subgraph of each $\mathcal{P}_k$ is $K_4 \oplus \tfrac{n-1}{2}Q_3$.

In each $\mathcal{P}_k$, let $u_k$ and $v_k$ be the two cells of the $K_4$ component in the $A_{\text{top-left}}$ and $A_{\text{bottom-right}}$ blocks, respectively. They have the same symbol. Deleting $u_k$ and $v_k$ leaves a single edge (the other two cells of the $K_4$) together with copies of $Q_3$. This graph is bipartite, so $\mathcal{P}_k \setminus \{u_k,v_k\}$ can be colored with two colors. Using a different pair of colors for each $k$, we color every cell outside $W = \{u_k, v_k : 0 \le k \le n-1\}$ with $2n$ colors.

We color $W$ with two new colors, red and blue. Let $U = \{u_k\}$ and $V = \{v_k\}$ be the cells of $W$ in the top-left and bottom-right blocks, respectively. Two observations describe how the cells of $W$ are adjacent.
\begin{itemize}
    \item The cell $v_k$ sits in exactly the same relative position as its counterpart $u_k$ (its row and column are those of $u_k$ shifted by $n$). So $v_k$ and $v_l$ share a row or column exactly when $u_k$ and $u_l$ do.

\item By the proof of \cref{lem:graph-dec}, $u_k$ and $v_k$ both carry the symbol $R_{2i^*_k + a_k} = R_{b_k - 1}$, since $2i^*_k \equiv b_k - a_k - 1 \pmod n$. As $(b_k)$ is a permutation, no two cells of $W$ from different $2$-plexes share a symbol. Moreover, a cell of $U$ and a cell of $V$ never share a row or column, so $u_k$ is adjacent to $v_l$ only when $l = k$.
\end{itemize}

Consider the graph induced on $U$. By \cref{lem2-plex}, every row and column contains at most two cells of $U$, so each cell has at most one neighbor in its row and at most one in its column. Thus the graph has maximum degree at most $2$ and is a disjoint union of paths and cycles. At each vertex of a cycle, one edge is a row edge and the other a column edge, so the edges alternate and the cycle has even length. Hence the graph induced on $U$ is bipartite, and we color it properly with red and blue. For $n=9$ this graph is a disjoint union of paths (see \cref{fig:d9_UV}).
Finally, give each $v_k$ the color opposite to that of $u_k$. Thus, $W$ is properly colored.

Altogether we used $2n+2$ colors, so $\chi(L_{D_n}) = 2n+2$.
\end{proof}

\begin{figure}[h]
\centering
\tiny
\setlength{\tabcolsep}{2.5pt}
\renewcommand{\arraystretch}{1.1}
\resizebox{\textwidth}{!}{%
\begin{tabular}{c|ccccccccc|ccccccccc}
$\cdot$ & $R_{0}$ & $R_{1}$ & $R_{2}$ & $R_{3}$ & $R_{4}$ & $R_{5}$ & $R_{6}$ & $R_{7}$ & $R_{8}$ & $S_{0}$ & $S_{8}$ & $S_{7}$ & $S_{6}$ & $S_{5}$ & $S_{4}$ & $S_{3}$ & $S_{2}$ & $S_{1}$ \\
\hline
$R_{0}$ & $R_{0}$ & $R_{1}$ & $R_{2}$ & $R_{3}$ & $R_{4}$ & $R_{5}$ & $R_{6}$ & $R_{7}$ & $R_{8}$ & $S_{0}$ & $S_{8}$ & $S_{7}$ & $S_{6}$ & $S_{5}$ & $S_{4}$ & $S_{3}$ & $S_{2}$ & $S_{1}$ \\
$R_{1}$ & $R_{1}$ & $R_{2}$ & $R_{3}$ & \cellcolor{softred}$\mathbf{R_{4}}$ & $R_{5}$ & $R_{6}$ & \cellcolor{skyblue}$\mathbf{R_{7}}$ & $R_{8}$ & $R_{0}$ & $S_{1}$ & $S_{0}$ & $S_{8}$ & $S_{7}$ & $S_{6}$ & $S_{5}$ & $S_{4}$ & $S_{3}$ & $S_{2}$ \\
$R_{2}$ & $R_{2}$ & \cellcolor{softred}$\mathbf{R_{3}}$ & $R_{4}$ & $R_{5}$ & $R_{6}$ & $R_{7}$ & $R_{8}$ & $R_{0}$ & $R_{1}$ & $S_{2}$ & $S_{1}$ & $S_{0}$ & $S_{8}$ & $S_{7}$ & $S_{6}$ & $S_{5}$ & $S_{4}$ & $S_{3}$ \\
$R_{3}$ & $R_{3}$ & $R_{4}$ & $R_{5}$ & $R_{6}$ & $R_{7}$ & $R_{8}$ & $R_{0}$ & $R_{1}$ & $R_{2}$ & $S_{3}$ & $S_{2}$ & $S_{1}$ & $S_{0}$ & $S_{8}$ & $S_{7}$ & $S_{6}$ & $S_{5}$ & $S_{4}$ \\
$R_{4}$ & $R_{4}$ & \cellcolor{skyblue}$\mathbf{R_{5}}$ & $R_{6}$ & $R_{7}$ & $R_{8}$ & $R_{0}$ & $R_{1}$ & \cellcolor{softred}$\mathbf{R_{2}}$ & $R_{3}$ & $S_{4}$ & $S_{3}$ & $S_{2}$ & $S_{1}$ & $S_{0}$ & $S_{8}$ & $S_{7}$ & $S_{6}$ & $S_{5}$ \\
$R_{5}$ & $R_{5}$ & $R_{6}$ & $R_{7}$ & \cellcolor{skyblue}$\mathbf{R_{8}}$ & $R_{0}$ & \cellcolor{softred}$\mathbf{R_{1}}$ & $R_{2}$ & $R_{3}$ & $R_{4}$ & $S_{5}$ & $S_{4}$ & $S_{3}$ & $S_{2}$ & $S_{1}$ & $S_{0}$ & $S_{8}$ & $S_{7}$ & $S_{6}$ \\
$R_{6}$ & $R_{6}$ & $R_{7}$ & $R_{8}$ & $R_{0}$ & $R_{1}$ & $R_{2}$ & $R_{3}$ & $R_{4}$ & $R_{5}$ & $S_{6}$ & $S_{5}$ & $S_{4}$ & $S_{3}$ & $S_{2}$ & $S_{1}$ & $S_{0}$ & $S_{8}$ & $S_{7}$ \\
$R_{7}$ & $R_{7}$ & $R_{8}$ & \cellcolor{softred}$\mathbf{R_{0}}$ & $R_{1}$ & $R_{2}$ & $R_{3}$ & $R_{4}$ & $R_{5}$ & \cellcolor{skyblue}$\mathbf{R_{6}}$ & $S_{7}$ & $S_{6}$ & $S_{5}$ & $S_{4}$ & $S_{3}$ & $S_{2}$ & $S_{1}$ & $S_{0}$ & $S_{8}$ \\
$R_{8}$ & $R_{8}$ & $R_{0}$ & $R_{1}$ & $R_{2}$ & $R_{3}$ & $R_{4}$ & $R_{5}$ & $R_{6}$ & $R_{7}$ & $S_{8}$ & $S_{7}$ & $S_{6}$ & $S_{5}$ & $S_{4}$ & $S_{3}$ & $S_{2}$ & $S_{1}$ & $S_{0}$ \\
\hline
$S_{0}$ & $S_{0}$ & $S_{8}$ & $S_{7}$ & $S_{6}$ & $S_{5}$ & $S_{4}$ & $S_{3}$ & $S_{2}$ & $S_{1}$ & $R_{0}$ & $R_{1}$ & $R_{2}$ & $R_{3}$ & $R_{4}$ & $R_{5}$ & $R_{6}$ & $R_{7}$ & $R_{8}$ \\
$S_{1}$ & $S_{1}$ & $S_{0}$ & $S_{8}$ & $S_{7}$ & $S_{6}$ & $S_{5}$ & $S_{4}$ & $S_{3}$ & $S_{2}$ & $R_{1}$ & $R_{2}$ & $R_{3}$ & \cellcolor{skyblue}$\mathbf{R_{4}}$ & $R_{5}$ & $R_{6}$ & \cellcolor{softred}$\mathbf{R_{7}}$ & $R_{8}$ & $R_{0}$ \\
$S_{2}$ & $S_{2}$ & $S_{1}$ & $S_{0}$ & $S_{8}$ & $S_{7}$ & $S_{6}$ & $S_{5}$ & $S_{4}$ & $S_{3}$ & $R_{2}$ & \cellcolor{skyblue}$\mathbf{R_{3}}$ & $R_{4}$ & $R_{5}$ & $R_{6}$ & $R_{7}$ & $R_{8}$ & $R_{0}$ & $R_{1}$ \\
$S_{3}$ & $S_{3}$ & $S_{2}$ & $S_{1}$ & $S_{0}$ & $S_{8}$ & $S_{7}$ & $S_{6}$ & $S_{5}$ & $S_{4}$ & $R_{3}$ & $R_{4}$ & $R_{5}$ & $R_{6}$ & $R_{7}$ & $R_{8}$ & $R_{0}$ & $R_{1}$ & $R_{2}$ \\
$S_{4}$ & $S_{4}$ & $S_{3}$ & $S_{2}$ & $S_{1}$ & $S_{0}$ & $S_{8}$ & $S_{7}$ & $S_{6}$ & $S_{5}$ & $R_{4}$ & \cellcolor{softred}$\mathbf{R_{5}}$ & $R_{6}$ & $R_{7}$ & $R_{8}$ & $R_{0}$ & $R_{1}$ & \cellcolor{skyblue}$\mathbf{R_{2}}$ & $R_{3}$ \\
$S_{5}$ & $S_{5}$ & $S_{4}$ & $S_{3}$ & $S_{2}$ & $S_{1}$ & $S_{0}$ & $S_{8}$ & $S_{7}$ & $S_{6}$ & $R_{5}$ & $R_{6}$ & $R_{7}$ & \cellcolor{softred}$\mathbf{R_{8}}$ & $R_{0}$ & \cellcolor{skyblue}$\mathbf{R_{1}}$ & $R_{2}$ & $R_{3}$ & $R_{4}$ \\
$S_{6}$ & $S_{6}$ & $S_{5}$ & $S_{4}$ & $S_{3}$ & $S_{2}$ & $S_{1}$ & $S_{0}$ & $S_{8}$ & $S_{7}$ & $R_{6}$ & $R_{7}$ & $R_{8}$ & $R_{0}$ & $R_{1}$ & $R_{2}$ & $R_{3}$ & $R_{4}$ & $R_{5}$ \\
$S_{7}$ & $S_{7}$ & $S_{6}$ & $S_{5}$ & $S_{4}$ & $S_{3}$ & $S_{2}$ & $S_{1}$ & $S_{0}$ & $S_{8}$ & $R_{7}$ & $R_{8}$ & \cellcolor{skyblue}$\mathbf{R_{0}}$ & $R_{1}$ & $R_{2}$ & $R_{3}$ & $R_{4}$ & $R_{5}$ & \cellcolor{softred}$\mathbf{R_{6}}$ \\
$S_{8}$ & $S_{8}$ & $S_{7}$ & $S_{6}$ & $S_{5}$ & $S_{4}$ & $S_{3}$ & $S_{2}$ & $S_{1}$ & $S_{0}$ & $R_{8}$ & $R_{0}$ & $R_{1}$ & $R_{2}$ & $R_{3}$ & $R_{4}$ & $R_{5}$ & $R_{6}$ & $R_{7}$ \\
\end{tabular}}

\vspace{10pt}
\begin{flushleft}
\hspace{6em}%
\begin{tikzpicture}
  \node[uvtx,fill=softred] (n0) at (0,0)   {$R_1$};
  \node[uvtx,fill=skyblue] (n2) at (1.8,0) {$R_8$};
  \node[uvtx,fill=softred] (n7) at (3.6,0) {$R_4$};
  \node[uvtx,fill=skyblue] (n4) at (5.4,0) {$R_7$};
  \node[uvtx,fill=softred] (n1) at (0,-1.4)   {$R_3$};
  \node[uvtx,fill=skyblue] (n3) at (1.8,-1.4) {$R_5$};
  \node[uvtx,fill=softred] (n6) at (3.6,-1.4) {$R_2$};
  \node[uvtx,fill=softred] (n5) at (0,-2.8)   {$R_0$};
  \node[uvtx,fill=skyblue] (n8) at (1.8,-2.8) {$R_6$};
  \begin{scope}[on background layer]
    \draw[rowe] (n0)--(n2) (n7)--(n4) (n3)--(n6) (n5)--(n8);
    \draw[cole] (n2)--(n7) (n1)--(n3);
  \end{scope}
\end{tikzpicture}

\vspace{6pt}
\hspace{6em}%
{\scriptsize
\tikz[baseline=-0.5ex]\draw[line width=1.5pt] (0,0)--(1em,0); same row\quad
\tikz[baseline=-0.5ex]\draw[line width=1.5pt,dash pattern=on 5pt off 3pt] (0,0)--(1.4em,0); same column}
\end{flushleft}

\caption{The Cayley table of $D_9$ (with $(a_k, b_k)$ as in \cref{tab:d9_2plexes}), with the cells of $U$ and $V$ in bold and the graph induced on $U$ shown below the table.}
\label{fig:d9_UV}
\end{figure}
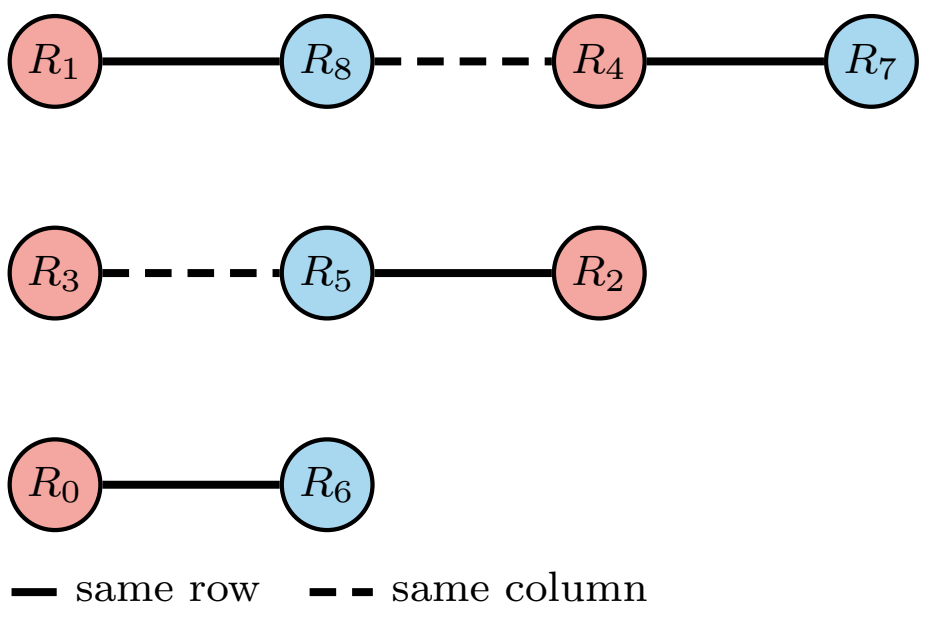


\section{Conclusion and future work} \label{future}
\cref{thm3} gives a fast and strong tool to determine if a finite group $G$ of order $n$ has a complete mapping and equivalently $\chi(L_G) = n$. In addition, \cref{thm8} guarantees that $\chi(L_G) \neq n + 1$ for any finite group $G$ of order $n$. Thus, for any finite group $G$ of order $n$ without a complete mapping $\chi(L_G) \geq n+2$.  \cref{thmabelian} shows that for any finite commutative group of order $n$ with a cyclic Sylow 2-subgroup, the chromatic number is $n + 2$. We prove that $\chi(L_{D_n}) = 2n + 2$ for an odd $n \geq 3$, which is the main result of this paper. Studying the chromatic number of $L_{D_{2n+1} \times \mathbb{Z}_{2m+1}}$ and $L_{Q_{2n+1}}$ has already begun and we have strong evidence that $\chi(L_{Q_{2n + 1}}) = 8n + 6$. Preparing a table of finite groups $G$ with orders less than $30$ and given chromatic numbers is in progress. In this table, there are $9$ groups without complete mapping, for which we did not prove that the chromatic number meets the lower bound. Proving \cref{conj1} for Latin squares of groups of order less than $30$ is included in our future work.

\section{Acknowledgments}
The first author is grateful to his MSc supervisor, the late Professor E.~S.~Mahmoodian of Sharif University of Technology, for introducing him to this problem during his MSc studies.

\bibliography{Ref}

\end{document}